# The Williams Conjecture is false for irreducible subshifts

By K. H. Kim and F. W. Roush*

**Abstract**

We show the Williams Conjecture is false for irreducible shifts of finite type by examining relative sign-gyration numbers of conjugacies between shifts with no points of period one or two.

## 1. Introduction

The shifts of finite type (SFTs) are topological dynamical systems which are the fundamental building blocks of symbolic dynamics, with applications to diverse topics such as smooth hyperbolic systems [1], coding [19], $C^*$-algebras [6], [9] and matrix theory [2]. The classification of SFTs (up to topological conjugacy) is a central open problem of symbolic dynamics. This problem has been dominated for over two decades by the Williams Conjecture that shift equivalence classifies SFTs. In [13], we gave a counterexample in the reducible case, but the conjecture remained open in the most important case, the irreducible case. In this paper we give the counterexample in the irreducible case announced in [14].

Both counterexamples grow out of the Factorization Theorem of [16], which related two fundamental representations of the automorphism group of a shift of finite type developed in the work of several authors. The Factorization Theorem states that the SGCC (Sign-Gyration-Compatability-Condition) representation of the automorphism group of an SFT factors through its dimension representation, and by certain explicit formulas. The proof of this theorem allows one to define a certain relative cohomology class whose nontriviality would give a counterexample to Williams' conjecture in the irreducible case (Section 7). In this paper we produce an example realizing this scheme, and in the course of this a simplification of the original proof of the Factorization Theorem (Section 6).

*The authors were partially supported by NSF Grants DMS 9024813 and DMS 9405004.

1991 *Mathematics Subject Classification.* Primary 58F03, 54H20.

To make the paper readily checkable, accessible and reasonably self-contained, and to exhibit the ultimate simplicity of the obstruction to Williams' conjecture, we provide sections of background and explanation through the developments of [16]. Although it was difficult for us to find any example at all to demonstrate the nontriviality of the obstruction, the example makes it clear that the obstruction is meaningful much more generally. The setting for the obstruction, which at present is precisely formulated but incompletely understood, allows optimism for further progress on the classification problem, which has been stuck for so long at the level of Williams' conjecture.

## 2. General background

A square nonnegative integral matrix $A$ presents an *edge shift of finite type* $\sigma_A$ as follows. Let $A$ be the adjacency matrix of a directed graph $\mathcal{G}$: If $A$ is $n \times n$, then $\mathcal{G}$ has vertices $1, 2, \ldots, n$ and the number of edges from vertex $i$ to vertex $j$ is $A(i,j)$. Let $X_A$ be the set of doubly infinite sequences $x = \ldots x_{-1} x_0 x_1 \ldots$ such that for all $i$, $x_i$ is in the edge set $\mathcal{E}$ of $\mathcal{G}$, and the terminal vertex of $x_i$ equals the initial vertex of $x_{i+1}$. Let $X_A$ have the compact, metrizable, zero-dimensional topology which is its relative topology as a subset of $\mathcal{E}^{\mathbb{Z}}$. Then $\sigma_A$ is the shift homeomorphism from $X_A$ to $X_A$, defined by $(\sigma_A x)_i = x_{i+1}$. A *topological conjugacy*, or isomorphism, from $\sigma_A$ to $\sigma_B$ is a homeomorphism $h: X_A \to X_B$ such that $\sigma_A h = h \sigma_B$. (In this paper, we will read composition from left to right.) The group of automorphisms (self-conjugacies) of $\sigma_A$ is denoted $\text{Aut}(\sigma_A)$. For a thorough introduction to SFTs, see [19].

Let $\Lambda$ be a subset of a ring such that $\Lambda$ contains 0 and 1. Let $A$ and $B$ be square matrices over $\Lambda$. An elementary strong shift equivalence over $\Lambda$ from $A$ to $B$ is a pair $(R, S)$ of rectangular matrices over $\Lambda$ such that $A = RS$ and $B = SR$. When such a pair exists, the matrices $A$ and $B$ are defined to be lag-1 strong shift equivalent over $\Lambda$. *Strong shift equivalence* is the equivalence relation on matrices over $\Lambda$ which is the transitive closure of lag-1 strong shift equivalence.

Matrices $A$ and $B$ with entries from $\Lambda$ are *shift equivalent* over $\Lambda$ if there exist matrices $R$ and $S$ over $\Lambda$ and a positive integer $n$ such that the following equations hold:

$$RA = BR, \quad AS = SB, \quad RS = B^n, \quad SR = A^n.$$

The ideas of shift equivalence and strong shift equivalence were introduced by Williams [26]. He proved that matrices $A$ and $B$ are strong shift equivalent over $\mathbb{Z}^+$ if and only if the SFTs $\sigma_A$ and $\sigma_B$ are isomorphic. However, strong shift equivalence over $\mathbb{Z}^+$ remains to this day a somewhat mysterious equivalence relation; for example, it is not known if there is a decision procedure for determining whether given matrices (even $2 \times 2$ matrices) are strong shift equivalent over $\mathbb{Z}^+$. Williams proposed shift equivalence over $\mathbb{Z}^+$ as a much more manageable and understandable equivalence relation.

For any $\Lambda$, strong shift equivalence over $\Lambda$ implies shift equivalence over $\Lambda$. Shift equivalence over $\mathbb{Z}^+$ can be understood and manipulated with semigroup theory, linear algebra, algebraic number theory and algebraic group theory, and there is a decision procedure [10], [11] for determining whether two matrices are shift equivalent over $\mathbb{Z}^+$. The final link of Williams' work [26] was the theorem: If two matrices are shift equivalent over $\mathbb{Z}^+$, then they are strong shift equivalent over $\mathbb{Z}^+$. Unfortunately, Williams' proof was incorrect; in the erratum, this theorem became Williams' conjecture [26], which has framed the study of the classification problem in the intervening years.

A matrix $A$ is *irreducible* if it is square and nonnegative and for every entry $(i,j)$ there exists $k > 0$ such that $A^k(i,j) > 0$. A matrix $A$ is *primitive* if there exists $k > 0$ such that every entry of $A^k$ is strictly positive. The topologically mixing SFTs are those which are isomorphic to $\sigma_A$ for primitive $A$. The SFTs with a dense forward orbit (the irreducible SFTs) are those which are isomorphic to $\sigma_A$ for irreducible $A$. Any SFT is the disjoint union of finitely many irreducible SFTs, which support all the ergodic theory and recurrent dynamics, and a wandering set of connecting orbits. The topological dynamics and classification of irreducible SFTs easily reduce to the mixing case. Thus Williams' conjecture is of greatest interest for irreducible matrices, especially primitive matrices.

Two facts (applied to $\Lambda = \mathbb{Z}$) greatly simplify the study of Williams' conjecture:

1. For $\Lambda$ a principal ideal domain: Shift equivalence over $\Lambda$ implies strong shift equivalence over $\Lambda$ [27], [3].

2. For $\Lambda$ a subring of the reals: Primitive matrices are shift equivalent over $\Lambda$ if and only if they are shift equivalent over $\Lambda^+$ ($:= \Lambda \cap [0,\infty)$) [21].

Thus in the primitive case, Williams' Conjecture can be reformulated as the statement: If two primitive matrices are strong shift equivalent over $\mathbb{Z}$, then they are strong shift equivalent over $\mathbb{Z}^+$. We will produce a counterexample to this formulation.

## 3. The $RS(\Lambda)$ complex

Again, $\Lambda$ is a subset of a ring such that $\Lambda$ contains 0 and 1. It is natural to view a strong shift equivalence over $\Lambda$ as a path. Wagoner introduced and developed the space of strong shift equivalences over $\Lambda$ as an algebraic topology framework for this view.

*Definition* 3.1 (1.4 of [23]; [24], [25]). The space $RS(\Lambda)$ of strong shift equivalences over $\Lambda$ is the geometric realization of the simplicial set where an $n$-simplex consists of the following data:

(a) an $(n+1)$-tuple $\langle A_0, \ldots, A_n \rangle$ of square matrices over $\Lambda$ and

(b) for each $i < j$ a strong shift equivalence $(R_{ij}, S_{ji})$ over $\Lambda$ from $A_i$ to $A_j$ such that the $RS$ Triangle Identities hold for $i < j < k$; that is,

$$R_{ij}R_{jk} = R_{ik}, \qquad R_{jk}S_{ki} = S_{ji}, \qquad S_{ki}R_{ij} = S_{kj}.$$

The face operators are the usual forgetful ones and the degeneracies insert the strong shift equivalence $(\mathrm{Id}, A_i)$ from $A_i$ to itself. The space $RS(\Lambda)$ is a naturally oriented CW complex. Abusing notation, we will refer, for example, to $A_i$ as a vertex, to $(R, S)$ as an edge and to $[(R_1, S_1), (R_2, S_2), (R_3, S_3)]$ as a triangle (for which the corresponding $\langle A_0, \ldots, A_n \rangle$ is $\langle R_1S_1, R_2S_2, R_3S_3 \rangle$).

To understand the genesis of the Triangle Identities, we must recall how an elementary strong shift equivalence $A = RS$, $B = SR$ over $\mathbb{Z}^+$ determines an *elementary conjugacy* $c$ from $\sigma_A$ and $\sigma_B$. View $A$ and $B$ as adjacency matrices of directed graphs with disjoint vertex sets, and view $R$ and $S$ as adjacency matrices for sets of edges between these vertex sets. According to the defining equations, for any vertices $i, j$ the number of $A$-edges $a$ from $i$ to $j$ equals the number of (two-edge) $RS$ paths $rs$ from $i$ to $j$ (here the initial vertex of the edge $r$ is $i$, the terminal vertex of the edge $s$ is $j$, and the terminal vertex of $r$ equals the initial vertex of $s$). So we may chose a bijection $\alpha$ of $A$-edges and $RS$-paths respecting initial and terminal vertices. Similarly we may choose a bijection $\beta$ of $B$-edges and $SR$-paths respecting initial and terminal vertices. Now the conjugacy $c$ is defined by the composition

$$\begin{aligned} x = \ldots x_{-1}x_0x_1\ldots &\mapsto \ldots (r_{-1}s_{-1})(r_0s_0)(r_1s_1)\ldots \\ &\mapsto \ldots (s_{-1}r_0)(s_0r_1)(s_1r_2)\ldots \\ &\mapsto \ldots y_0y_1y_2\ldots = y \end{aligned}$$

where the first and third bijections are defined by coordinatewise application of the bijections $\alpha$ and $\beta$, and the middle bijection just shifts parentheses.

If $\pi$ is a permutation of edges respecting initial and terminal vertices, then the map $x \mapsto x'$ given by

$$\begin{aligned} &\ldots x_{-1}x_0x_1\ldots \\ \mapsto &\ldots (\pi(x_{-1}))(\pi(x_0))(\pi(x_1))\ldots \\ = &\ldots x'_{-1}x'_0x'_1\ldots \end{aligned}$$

determines a self conjugacy (automorphism). Such an automorphism is an elementary *simple* automorphism of Nasu [20]. An automorphism of $\sigma_A$ is *simple* if it is a composition of automorphisms of the form $\phi c \phi^{-1}$ where $c$ is elementary simple and $\phi$ is a conjugacy. Note, the conjugacy $c$ of the previous paragraph is uniquely determined by $(R, S)$ up to composition by simple automorphisms.

A triangle in $RS(\mathbb{Z}^+)$ has a crucial property: there are associated elementary conjugacies $c_i = c(R_i, S_i)$ such that $c_1 c_2 = c_3$ modulo composition by simple automorphisms [25]. The "modulo" qualification in this statement is necessary for a single reason: if two of the three edges are the same, then it may be necessary to choose different bijections of edges in the definitions of the associated conjugacies in order to get a commuting diagram of conjugacies. For example, if $(R_1, S_1) = (R_3, S_3)$ and $(R_2, S_2)$ forces $c_2$ to be a nontrivial automorphism, then the choices of bijections for $c_1$ and $c_3$ must be different to permit $c_1 c_2 = c_3$.

The Triangle Identities were designed so that two paths from $A$ to $B$ in $RS(\{0, 1\})$ give rise to the same conjugacy if and only if the paths are homotopic in $RS(\{0, 1\})$ [24]. In particular, $\pi_1(RS(\{0, 1\}), A)$, the fundamental group of a component of $RS(\{0, 1\})$ with basepoint $A$, is isomorphic to $\text{Aut}(\sigma_A)$. An important consequence is that two paths from $A$ to $B$ are homotopic in $RS(\mathbb{Z}^+)$ if and only if they give rise to conjugacies which are equal modulo composition with a simple automorphism [25]. In particular, $\pi_1(RS(\mathbb{Z}^+), A) \cong \text{Aut}(\sigma_A)/\text{Simp}(\sigma_A)$, where $\text{Simp}(\sigma_A)$ is the group of simple automorphisms of $\sigma_A$. Thus one has an algebraic topology setting for studying the automorphisms of $\sigma_A$, which has proved very useful. For purposes of our counterexample to Williams' conjecture, we do not need this deeper theory, only more direct consequences of the triangle property of the preceding paragraph.

## 4. Direct limit modules and $RS(\Lambda)$

Suppose $\Lambda$ is a ring with unity, acting from the left on $\Lambda^n$, the free $\Lambda$-module of rank $n$ with elements given as row vectors. Let $A$ be an $n \times n$ matrix over $\Lambda$. Let $G_A$ denote the direct limit $\Lambda$-module obtained by the action of $A$ from the right on $\Lambda^n$. So, an element of $G_A$ is an equivalence class $[(v, i)]$, where $v \in \mathbb{Z}^n, i \in \mathbb{N}$ and $[(v, i)] = [(w, j)]$ if and only if $vA^{j+k} = wA^{i+k}$ for some $k > 0$. The rule $[(v, i)] \mapsto [(vA, i)]$ defines an automorphism $s_{A/\Lambda}$ of the module $G_A$. We let $\text{Aut}(s_{A/\Lambda})$ denote the group of $\Lambda$-module automorphisms of $G_A$ which commute with $s_{A/\Lambda}$. If $\Lambda = \mathbb{Z}$, then the module $G_A$ is just an abelian group, and we may suppress $\Lambda$ from the notation. It is well-known that square matrices over $\Lambda$ define isomorphic $\Lambda$-modules as above if and only if the matrices are shift equivalent over $\Lambda$ ([18], [3]).

If $A$ is nonsingular and $\Lambda$ is an integral domain, then $\text{Aut}(s_{A/\Lambda})$ is simply the group of nonsingular matrices $M$ over the field of fractions of $\Lambda$ such that $AM = MA$ and for all large $k$, $MA^k$ has all entries in $\Lambda$. Here the action of a matrix $R$ is by $\widehat{R} : [(v, i)] \mapsto [(vR, i)]$.

An elementary strong shift equivalence $(R, S)$ over $\Lambda$ from $A = RS$ to $B = SR$ induces the isomorphism of $\Lambda$-modules $\widehat{R} : [(v, i)] \mapsto [(vR, i)]$. Now let $\mathcal{P}$ denote a path $(R_1, S_1)^{\varepsilon(1)} (R_2, S_2)^{\varepsilon(2)} \ldots (R_k, S_k)^{\varepsilon(k)}$ of edges from $A$ to $B$ in $RS(\Lambda)$, where $\varepsilon(i)$ is 1 or $-1$ depending on whether the edge $(R_i, S_i)$ is traversed with positive or negative orientation. Any path from $A$ to $B$ is homotopic to some such path. We associate to $\mathcal{P}$ the module isomorphism $\widehat{\mathcal{P}} : G_A \to G_B$ which is the composition $(\widehat{R}_1)^{\varepsilon(1)} \ldots (\widehat{R}_k)^{\varepsilon(k)}$.

Finally, suppose that $\Lambda$ is a principal ideal domain. For this case we can summarize the significance of the Triangle Identities in $RS(\Lambda)$ as derived by Wagoner [24]. The automorphism $\widehat{\mathcal{P}}$ described above depends only on the homotopy class of the path $\mathcal{P}$ from $A$ to $B$. Also, the map $\mathcal{P} \mapsto \widehat{\mathcal{P}}$ induces an isomorphism from $\pi_1(RS(\Lambda), A) \to \text{Aut}(s_{A/\Lambda})$.

For a matrix $A$ over $\mathbb{Z}^+$, it follows that the inclusion map $RS(\mathbb{Z}^+) \hookrightarrow RS(\mathbb{Z})$ induces a homomorphism $\text{Aut}(\sigma_A)/\text{Simp}(\sigma_A) \to \text{Aut}(s_A)$. Precomposition with the natural map $\text{Aut}(\sigma_A) \to \text{Aut}(\sigma_A)/\text{Simp}(\sigma_A)$ gives a presentation of Krieger's *dimension representation* $\rho : \text{Aut}(\sigma_A) \to \text{Aut}(s_A)$ [4], [5].

## 5. Periodic points and relative sign-gyration numbers

Suppose $\phi$ is a conjugacy from $\sigma_A$ to $\sigma_B$. For a given positive integer $m$, let $P_m^o(\sigma_A)$ denote the points in $\sigma_A$-orbits of cardinality $m$. Let $(x_1, \ldots, x_k)$ be a basis for $P_m^o(\sigma_A)$: that is, a tuple of distinct points $x_i$ such that each $\sigma_A$-orbit of cardinality $m$ contains exactly one of the points $x_i$. Similarly let $(y_1, \ldots, y_k)$ be a basis for $P_m^o(\sigma_B)$. Then there is a permutation $\pi$ of $\{1, \ldots, k\}$ and a $k$-tuple of integers $(n(1), \ldots, n(k))$ such that $\phi(x_i) = (\sigma_B)^{n(i)} y_{\pi(i)}$. We define the *orbit sign number* $\text{OS}_n(\phi)$ in $\mathbb{Z}/2$ to be 0 if the permutation $\pi$ has even sign and 1 otherwise. We define the *gyration number* in $\mathbb{Z}/m$ to be $\text{GY}_m(\phi) = \sum_i n(i)$. Finally we define in $\mathbb{Z}/m$ the *sign-gyration* number

$$\text{SGCC}_m(\phi) = \text{GY}_m(\phi) + (m/2) \sum_i \text{OS}_i(\phi)$$

where the last sum is over the positive integers $i < m$ such that $m/i$ is a power of 2 (and an empty sum is 0). "SGCC" is as defined in the introduction [4].

If $\phi$ is an automorphism (i.e., $A = B$) and the two bases $(x_1, \ldots, x_k)$, $(y_1, \ldots, y_k)$ are the same, then $\phi \mapsto \text{OS}_m(\phi)$, $\phi \mapsto \text{GY}_m(\phi)$ and $\phi \mapsto \text{SGCC}_m(\phi)$ define group homomorphisms from $\text{Aut}(\sigma_A)$ to $\mathbb{Z}/2$ and $\mathbb{Z}/m$ which do not depend on the particular choice of base. We let GY, OS and SGCC denote the product maps, e.g. $\text{GY}(\phi) = \prod_m \text{GY}_m(\phi) \in \prod_n \mathbb{Z}/m$. A key fact for us is that SGCC vanishes on simple automorphisms [20].

To apply these ideas in $RS(\mathbb{Z}^+)$, at each vertex $A$ we make a definite choice of basis for each $P_n^o(\sigma_A)$, and at each edge $(R,S)$ we make a definite choice of elementary conjugacy $c(R,S)$. This is done by certain natural lexicographic rules, chosen in [16]. For each edge $(R,S)$, these choices determine $\mathrm{SGCC}_m(c(R,S))$, which we abbreviate as $\mathrm{SGCC}_m(R,S)$. To a path $\mathcal{P} = (R_1,S_1)^{\varepsilon(1)}(R_2,S_2)^{\varepsilon(2)}\ldots(R_k,S_k)^{\varepsilon(k)}$ from $A$ to $B$ in $RS(\mathbb{Z}^+)$ we associate the *relative sign-gyration number*

$$\mathrm{SGCC}_m(\mathcal{P}) = \sum_i \varepsilon(i)\mathrm{SGCC}_m(R_i,S_i) \ .$$

From the explicit combinatorics, one derives [16] for each $m$ a function $\mathrm{sgc}_m$ from edges $(R,S)$ of $RS(\mathbb{Z}^+)$ into $\mathbb{Z}/m$ such that $\mathrm{sgc}_m(R,S) = \mathrm{SGCC}_m(R,S)$, giving

$$\mathrm{SGCC}_m(\mathcal{P}) = \sum_i \varepsilon(i)\mathrm{sgc}_m(R_i,S_i) \ .$$

The function $\mathrm{sgc}_m$ is a polynomial function of the entries of $R$ and $S$ with rational coefficients whose denominators divide $(2m)!$. In particular, $\mathrm{sgc}_m(R,S)$ only depends on the the values of $R$ and $S$ modulo $m(2m)!$ [16]. As $m$ increases, $\mathrm{sgc}_m$ becomes extremely complicated. Our counterexample will use only $\mathrm{sgc}_2$, which is given by the formula

$$\mathrm{sgc}_2(R,S) = \sum_{\substack{i<j \\ k>l}} R_{ik}S_{ki}R_{jl}S_{lj} + \sum_{\substack{i<j \\ k\geq l}} R_{ik}S_{kj}R_{jl}S_{li} + \sum_{i,k} \frac{R_{ik}(R_{ik}-1)}{2}S_{ki}^2 \ .$$

Note that $\mathrm{sgc}_2(R,S)$ only depends on the values of $R$ and $S$ modulo 4.

Now suppose $c_1, c_2, c_3$ are the elementary conjugacies we have associated to the sides of an $RS(\mathbb{Z}^+)$ triangle. If $c_1c_2 = c_3$, then a simple computation [16] shows

$$\mathrm{SGCC}(R_1,S_1) + \mathrm{SGCC}(R_2,S_2) = \mathrm{SGCC}(R_3,S_3) \ .$$

Because SGCC vanishes on simple automorphism and $c_1c_2 = c_3$ modulo simple automorphisms, this addition formula holds for all $RS(\mathbb{Z}^+)$ triangles. That is, SGCC vanishes around $RS(\mathbb{Z}^+)$ triangles.

## 6. The obstruction

In this section, $m$ is a positive integer and $L_m = L$ denotes a positive integer such that $L/m$ is an integer which is divisible by the denominators of the coefficients in the polynomial formula $\mathrm{sgc}_m$. If $m = 2$, then we can take $L = 4$; for any $m$, we could use $L = m(2m)!$. We define $\mathrm{sgc}_m$ from $RS(\Lambda)$-paths into $\mathbb{Z}/m$ by the same polynomial formulas as for $RS(\mathbb{Z}^+)$-paths, when the formulas make sense on $\Lambda$. The formulas do make sense if $\Lambda = \mathbb{Z}/L$. Likewise the formulas make sense if $\Lambda$ is a set of rational numbers whose denominators are relatively prime to $L$, if we view $p/q$ as $pq^{-1}$ (modulo $L$). For $\Lambda = \mathbb{Z}$ we define $\mathrm{sg}c = \prod_m \mathrm{sg}c_m$ as a homomorphism into $\prod_m \mathbb{Z}/m$.

SGCC is a function on paths in $RS(\mathbb{Z}^+)$ which carries information about the action of associated topological conjugacies on periodic points. The crucial issue is to relate this information to the algebra of matrices and eigenvalues surrounding the direct limit modules. To this end, we want to extend SGCC to a function on paths of edges in $RS(\mathbb{Z})$ which is compatible with the $RS(\mathbb{Z})$ Triangle Identities. We will check $\mathrm{sg}c$ vanishes around $RS(\mathbb{Z})$ triangles. Because $\mathrm{sg}c$ agrees with SGCC on $RS(\mathbb{Z}^+)$, we can then simply use $\mathrm{sg}c$ for the desired extension of SGCC. It will be computationally useful to prove the more general statement:

COCYCLE LEMMA 6.1.  *Suppose $\Lambda$ is* (i) $\mathbb{Z}/L$ *or* (ii) *a set of rationals whose denominators are relatively prime to $L$. Then $\mathrm{sg}c_m$ vanishes around $RS(\Lambda)$ triangles.*

*Proof.* Suppose that $T = [(R_1, S_1), (R_2, S_2), (R_3, S_3)]$ is an $RS(\Lambda)$ triangle. In the case $\Lambda = \mathbb{Z}/L$, choose positive integral matrices $R_1', R_2', S_3'$ congruent (modulo $L$) respectively to $R_1, R_2, S_3$. Now define matrices $R_3', S_1', S_2'$ from these three matrices with the Triangle Identities:

$$R_3' := R_1'R_2', \qquad S_1' := R_2'S_3', \qquad S_2' := S_3'R_1'.$$

This gives us an $RS(\mathbb{Z}^+)$ triangle $T' = [(R_1', S_1'), (R_2', S_2'), (R_3', S_3')]$. Then $\mathrm{sg}c_m$ vanishes around $T'$, because $\mathrm{sg}c_m = \mathrm{SGCC}_m$ around $T'$. Therefore $\mathrm{sg}c_m$ vanishes around $T$, because $\mathrm{sg}c_m$ only depends on the matrix entries modulo $L$.

In case (ii), $T$ modulo $L$ is an $RS(\mathbb{Z}/L)$ triangle, and we reduce to case (i). □

The Factorization Theorem (1.4 of [16]) stated that the SGCC homomorphism on $\mathrm{Aut}(\sigma_A)$ factored through the dimension representation $\rho$, by certain explicit formulas. We want to point out here that we have proved a version of the Factorization Theorem (which improves a little on the original, as explained in Remark 8.1). We present the dimension representation $\rho : \mathrm{Aut}(\sigma_A) \to \mathrm{Aut}(s_A)$ as $\rho : \mathrm{Aut}(\sigma_A) \to \pi_1(RS(\mathbb{Z}), A)$, as explained in Section 4. Now the homomorphisms $\mathrm{sg}c_m$ on paths of edges in $RS(\mathbb{Z})$ depend only on the homotopy class of a path and therefore induce a homomorphism from $\pi_1(RS(\mathbb{Z}), A)$ into $\prod_m \mathbb{Z}/m$. We regard this as a map (which for simplicity we also call $\mathrm{sg}c$) from $\mathrm{Aut}(s_A)$ into $\prod_m \mathbb{Z}/m$. Then we can summarize the improved version of the Factorization Theorem as follows.

FACTORIZATION THEOREM 6.2.  *If $A$ is a square matrix over $\mathbb{Z}^+$, then the map SGCC on $\mathrm{Aut}(\sigma_A)$ is the composition $\rho$ followed by $\mathrm{sg}c$.*

Now suppose the Cocycle Lemma applies to a PID $\Lambda$ and $\mathcal{P}$ is a strong shift equivalence over $\Lambda$ from $A$ to $B$. We view $\mathcal{P}$ as a path of edges in $RS(\Lambda)$. The relative sign-gyration number $\mathrm{sgc}_m(\mathcal{P})$ depends only on the homotopy class of $\mathcal{P}$ as a path in $RS(\Lambda)$. Identify $\pi_1(RS(\Lambda), A)$ with $\mathrm{Aut}(s_{A/\Lambda})$. Then we can associate to the pair $(A, B)$ the subset $\mathrm{rsgc}_{m,\Lambda}[A, B]$ of $\mathbb{Z}/m$ which is the ($m^{\mathrm{th}}$) *relative sign-gyration* coset

$$\mathrm{rsgc}_{m,\Lambda}[A, B] = \mathrm{sgc}_m(\mathcal{P}) + \mathrm{sgc}_m\left(\mathrm{Aut}(s_{A/\Lambda})\right).$$

This set depends only on the pair $(A, B)$ and we can now formulate the obstruction: if $\mathcal{P}'$ is any strong shift equivalence over $\Lambda$ from $A$ to $B$, then $\mathrm{sgc}_m(\mathcal{P}')$ must lie in $\mathrm{rsgc}_{m,\Lambda}[A, B]$.

When $\Lambda = \mathbb{Z}$, we may omit it from the notation as before. In the case $\Lambda = \mathbb{Z}$, we may combine the information over all $m$ and define the *relative sign-gyration coset* in $\mathrm{prod}_m \mathbb{Z}/m$ as

$$\mathrm{rsgc}[A, B] = \mathrm{sgc}(\mathcal{P}) + \mathrm{sgc}\left(\mathrm{Aut}(s_A)\right).$$

## 7. The counterexample

We will produce primitive matrices $A$ and $B$ satisfying the following conditions:

(1) $\mathrm{tr}(A) = \mathrm{tr}(A^2) = 0$,

(2) $\mathrm{sgc}_2$ vanishes on $\mathrm{Aut}(s_A)$ and

(3) there is a path in $RS(\mathbb{Z})$ from $A$ to $B$ with nonzero $\mathrm{sgc}_2$.

It follows from (2) and (3) that $\mathrm{rsgc}_2[A, B] = \{1\}$, so $\mathrm{sgc}_2(\mathcal{P}) = 1$ for every strong shift equivalence $\gamma$ from $A$ to $B$. If $\mathcal{P}$ were a strong shift equivalence over $\mathbb{Z}^+$ from $A$ to $B$, then $\mathrm{sgc}_2(\mathcal{P})$ would vanish, since SGCC = sgc over $\mathbb{Z}^+$ and the subshifts $\sigma_A, \sigma_B$ have no points of period 1 or 2. Therefore $A$ and $B$ are not strong shift equivalent over $\mathbb{Z}^+$.

We can return to the language of the introduction to describe this. Let $RS_m(\mathbb{Z}^+)$ denote the union of path components of $RS(\mathbb{Z}^+)$ containing those primitive $A$ having no points of period less than or equal to $m$. The map $\mathcal{P} \mapsto \mathrm{sgc}_m(\mathcal{P})$ gives a cohomology class in $H^1(RS(\mathbb{Z}), RS_m(\mathbb{Z}^+); \mathbb{Z}/m)$. If this class is nonzero on any $RS(\mathbb{Z})$ path between primitive $A$ and $B$ where $\mathrm{tr}(A^m) = 0$, then $A$ and $B$ are shift equivalent but not strong shift equivalent over $\mathbb{Z}^+$.

It remains to actually find matrices satisfying the conditions. This was difficult for us ([14], Section 5). Define

$$S = \begin{pmatrix} 2 & 2 & 2 & 1 & 3 & 0 & 0 \\ 1 & 2 & 2 & 1 & 3 & 0 & 0 \\ 1 & 1 & 2 & 1 & 3 & 0 & 0 \\ 1 & 1 & 1 & 1 & 3 & 0 & 0 \\ 0 & 0 & 0 & 0 & 0 & 0 & 1 \\ 4 & 5 & 6 & 3 & 10 & 0 & 0 \\ 4 & 5 & 6 & 3 & 0 & 1 & 0 \end{pmatrix}, \quad R = \begin{pmatrix} -1 & 0 & 1 & 1 & 0 & 0 & 0 \\ 1 & -1 & 0 & 0 & 0 & 0 & 0 \\ 0 & 1 & -1 & 0 & 0 & 0 & 0 \\ 0 & 0 & 1 & -1 & 0 & 0 & 0 \\ 0 & 0 & 0 & 0 & 1 & 0 & 0 \\ 0 & 0 & 0 & 0 & 0 & 1 & 0 \\ 0 & 0 & 0 & 0 & 0 & 0 & 1 \end{pmatrix},$$

$$A = \begin{pmatrix} 0 & 0 & 1 & 1 & 3 & 0 & 0 \\ 1 & 0 & 0 & 0 & 3 & 0 & 0 \\ 0 & 1 & 0 & 0 & 3 & 0 & 0 \\ 0 & 0 & 1 & 0 & 3 & 0 & 0 \\ 0 & 0 & 0 & 0 & 0 & 0 & 1 \\ 1 & 1 & 1 & 1 & 10 & 0 & 0 \\ 1 & 1 & 1 & 1 & 0 & 1 & 0 \end{pmatrix}, \quad B = \begin{pmatrix} 0 & 0 & 1 & 1 & 3 & 0 & 0 \\ 1 & 0 & 0 & 0 & 0 & 0 & 0 \\ 0 & 1 & 0 & 0 & 0 & 0 & 0 \\ 0 & 0 & 1 & 0 & 0 & 0 & 0 \\ 0 & 0 & 0 & 0 & 0 & 0 & 1 \\ 4 & 5 & 6 & 3 & 10 & 0 & 0 \\ 4 & 5 & 6 & 3 & 0 & 1 & 0 \end{pmatrix}.$$

The matrices $A$ and $B$ are primitive; $A = SR, B = RS$; and $\text{tr}(A) = \text{tr}(A^2) = 0$. The condition $\text{sg}c_2(R,S) \neq 0$ can be checked with a simple program (for example the Maple program provided in [14]). It remains to check that $\text{sg}c_2$ vanishes on $\text{Aut}(s_A)$. The determinant of $A$ is $-1$, so $\text{Aut}(s_A)$ is $\mathcal{C}(A)$, the centralizer of $A$ in $\text{GL}(7,\mathbb{Z})$. The characteristic polynomial of $A$ is $p(t) = 1 - 17t - 33t^2 - 28t^3 - 23t^4 + t^7$. Because $p(t)$ is irreducible, $\mathcal{C}(A) \subset \mathbb{Q}[A]$ and $\mathbb{Q}[A]$ is isomorphic to the algebraic number field $\mathbb{Q}[t]/p(t)$ under the isomorphism induced by $t \mapsto A$. Under this isomorphism, $\mathcal{C}(A)$ corresponds to a subgroup of the units group $\mathcal{U}$ of the algebraic integers of this field. Because $\mathbb{Q}(A)$ has a real embedding, the only torsion elements of $\mathcal{U}$ are 1 and $-1$. Because $\text{tr}(A) = 0$, the formula for $\text{sg}c_2$ shows $\text{sg}c_2(-I) = 0$.

The polynomial $p$ has three real and four complex roots, so the Dirichlet Units Theorem gives $\mathcal{U} \cong \mathbb{Z}^4 \oplus \mathbb{Z}/2$. A PARI computer calculation gives the following system of fundamental units for $\mathcal{U}$:

$$\begin{aligned} f_1 &= t, \\ f_2 &= \frac{1}{3}(38t^6 + 2t^5 - t^4 - 872t^3 - 1108t^2 - 1309t - 713), \\ f_3 &= \frac{1}{3}(842t^6 + 5072t^5 - 6847t^4 - 46061t^3 - 34930t^2 - 52216t + 2878), \\ f_4 &= \frac{1}{3}(4260971t^6 - 3124108t^5 + 2290532t^4 - 99681839t^3 - 46221667t^2 \\ &\quad - 106722952t + 5811547). \end{aligned}$$

For each of these polynomials $f_i$, we compute an edge from $A$ to $A$ in $RS(\mathbb{Z}[1/3])$ given by $(R,S) = (f_i(A), [f_i(A)]^{-1}A)$. Because $\mathcal{C}(A)$ corresponds to a proper subgroup of $\mathcal{U}$, the matrices $f_i(A)$ are not all integral. Nevertheless, $\mathcal{C}(A)$ is in the group generated under multiplication by the $f_i(A)$ and $-I$; an $RS(\mathbb{Z})$ triangle is an $RS(\mathbb{Z}[1/3])$ triangle; and $\mathrm{sg}c_2$ is well-defined on homotopy classes of paths in $RS(\mathbb{Z}[1/3])$. So to show that $\mathrm{sg}c_2$ vanishes on $\mathrm{Aut}(s_A)$, it suffices now to check that $\mathrm{sg}c_2$ vanishes on the four edges $(R,S) = (f_i(A), [f_i(A)]^{-1}A)$. For this we multiply the matrices by 9 to clear denominators ($9 \equiv 1 \bmod 4$) and apply our little program. That finishes the proof.

We remark that in checking, we can avoid dependence on PARI. We can variously check that a polynomial $f(t)$ indeed defines an algebraic unit, e.g. crudely by examining the characteristic polynomial of $f(A)$. Then we only need to verify that the system of four units together with $\{-1\}$ is a system of generators for $\mathcal{U}$ modulo squares. For this it suffices to define a homomorphism $\pi$ from $\mathcal{U}$ to $(\mathbb{Z}/2)^5$ which sends $(f_1(t), f_2(t), f_3(t), f_4(t), -1)$ to a basis. We choose $\pi$, the product of homomorphisms $\pi_1, \ldots, \pi_5$ into $\mathbb{Z}/2$. For each $\pi_i$, we choose an odd prime $m$ and an integer $t$ such that $p(t) \equiv 0 \pmod{m}$; we define $\pi_i(f) = 0$ if $f(t)$ is a quadratic residue $\pmod{m}$ and $\pi_i(f) = 1$ otherwise. We choose $(m_1, m_2, m_3, m_4, m_5) = (17, 17, 41, 11, 11)$ and $(t_1, t_2, t_3, t_4, t_5) = (2, 4, 3, 1, 10)$ and summarize the computations in the matrices below, where $M(i,j) = f_i(t_j)$ and $Q(i,j) = \pi_i(f_j)$. The matrix $Q$ is invertible (mod 2).

$$M = \begin{pmatrix} 2 & 5 & 9 & 10 & -1 \\ 4 & 9 & 5 & 8 & -1 \\ 3 & 26 & 2 & 36 & -1 \\ 1 & 10 & 4 & 8 & -1 \\ 10 & 7 & 9 & 6 & -1 \end{pmatrix}, \quad Q = \begin{pmatrix} 0 & 1 & 0 & 1 & 0 \\ 0 & 0 & 1 & 0 & 0 \\ 1 & 1 & 0 & 0 & 0 \\ 0 & 1 & 0 & 1 & 1 \\ 1 & 1 & 0 & 1 & 1 \end{pmatrix}.$$

## 8. Remarks

1. The Cocycle Lemma showed $\mathrm{sg}c_m(R,S)$ is an extension of $\mathrm{SGCC}_m(R,S)$ vanishing around $RS(\mathbb{Z})$ triangles. The localization and positivity arguments of [16] produced a different extension $\mathrm{sgcc}_m(R,S)$ of $\mathrm{SGCC}_m(R,S)$ vanishing around $RS(\mathbb{Z})$ triangles. The function $\mathrm{sgcc}_m(R,S)$ is defined when the characteristic polynomial of $A$ has a simple positive root strictly greater than the modulus of any other root, as it must when $A$ is shift equivalent over $\mathbb{Z}$ to a primitive matrix. In this case (see (3.3) of [16]),

$$\begin{aligned} \mathrm{sgcc}_m(R,S) &= \mathrm{sg}c_m(R,S) \quad \text{if } \mathrm{sgn}(R) > 0, \text{ and} \\ \mathrm{sgcc}_m(R,S) &= \mathrm{sg}c_m(-R,-S) \quad \text{if } \mathrm{sgn}(R) < 0 \end{aligned}$$

where $\mathrm{sgn}(R)$ denotes the sign of the number by which $R$ multiplies the eigenvector of the dominant root.

Thus the version of the Factorization Theorem proved in Section 6 differs from the original in two ways: $A$ is no longer required to be primitive, and in the primitive case the explicit homomorphism sg$c$ differs from the explicit homomorphism used in [16]. The proof differs in that the Cocycle Lemma is considerably simpler than the localization and positivity arguments of [16]. However, our counterexample would work just as well with the extension sg$cc_2$. The simplified extension with the Cocycle Lemma was pointed out to us by M. Boyle in the course of examining our original version.

2. J. Wagoner has pointed out to us that the functions sg$cc_2$ and sg$c_2$ do not agree in general, for example on the edge $([-1], [-3])$, but they do agree on components of $RS(\mathbb{Z})$ containing a vertex which is a primitive matrix with trace zero. To see this, note that the Cocycle Lemma implies

(i) sg$c_2(-I, -A) =$ sg$c_2(-I, -B)$ if $A$ and $B$ are SSE over $\mathbb{Z}$.

(ii) sg$c_2(-R, -S) =$ sg$c_2(R, S) +$ sg$c_2(-I, -A)$.

If A is nonnegative and tr$(A) = 0$, then the formula for sg$c_2$ shows sg$c_2(-I, -A) = 0$. Hence sg$c_2(-I, -M) = 0$ for any $M$ which is strong shift equivalent to $A$ over $\mathbb{Z}$, and the difference between sg$cc_2$ and sg$c_2$ disappears.

3. Computer calculations of sg$c_m$ rapidly approach impossibility as $m$ increases, even when one uses some other special formulas adapted to the case when $R, S$ are polynomials in a common variable.

4. When $A$ has entries in $\mathbb{Z}^+$, $G_A$ acquires a natural order structure making it a dimension group. (For this reason, the unordered group $G_A$ is sometimes referred to as the dimension group of $A$ or of $\sigma_A$ [19].) This also explains the terminology for Krieger's dimension representation, which is defined by way of a Grothendieck-style construction of a version of $G_A$ from certain compact sets [18], [4], [5]. In our paper the order structure plays no role.

5. The Factorization Theorem was originally proved as a key ingredient for understanding the action of Aut$(\sigma_A)$ on periodic points. The condition SGCC $= 0$ is the only obstruction to the action of Ker(Aut$(\sigma_A)$) on finite subsystems of $\sigma_A$ [17].

6. The counterexample [13] to the reducible Williams Conjecture involved a fundamentally different argument which relied on one consequence of the Factorization Theorem: the example (4.1) of [16] for which the dimension representation is not surjective. More generally, one ingredient which is necessary for the classification of reducible SFTs is the determination of the range of the dimension representation on irreducible SFTs [15].

7. Although we have produced a counterexample to the irreducible Williams Conjecture, it is in no way a repudiation of the ideas of shift equivalence and strong shift equivalence, which remain part of the foundation for work on the classification problem.

*Acknowledgment.* We have learned that Jack Wagoner developed this basic strategy for finding a counterexample and the theory of the invariant in unpublished work prior to our own (independent) discovery of these ideas. We thank Mike Boyle for the simplifying Cocycle Lemma and Jack Wagoner for the remarks clarifying the relation of sg$c$ and SGCC, and we thank both of them for substantial help with the exposition.

Mathematics Research Group, Alabama State University, Montgomery, AL.
and Korean Academy of Science and Technology
*E-mail address*: kkim@asu.alasu.edu
*E-mail address*: froush@asu.alasu.edu

References

[1] R. Bowen, *Equilibrium States and the Ergodic Theory of Anosov Diffeomorphisms*, Lecture Notes in Math. **470** (1975), Springer-Verlag, New York.
[2] M. Boyle, Symbolic dynamics and matrices, in *Combinatorial and Graph-Theoretical Problems in Linear Algebra*, I.M.A. Vol. Math. Appl. **50** (1993), 1–38, Springer-Verlag, New York.
[3] M. Boyle and D. Handelman, Algebraic shift equivalence and primitive matrices, Trans. AMS **336** (1993), 121–149.
[4] M. Boyle and W. Krieger, Periodic points and automorphisms of the shift, Trans. AMS **302** (1987), 125–149.
[5] M. Boyle, D. Lind, and D. Rudolph, The automorphism group of a shift of finite type, Trans. AMS **306** (1988), 71–114.
[6] J. Cuntz and W. Krieger, A class of C*-algebras and topological Markov chains, Invent. Math. **56** (1980), 251–268.
[7] U. Fiebig, Gyration numbers for involutions of subshifts of finite type I, Forum Math.4 (1992), 77-108 and II, Forum Math. **4** (1992), 183–211.
[8] Y. W. Ha, Conjugation and strong shift equivalence, Commun. Korean Math. Soc. **11** (1996), 191–199.
[9] D. Huang, Flow equivalence of reducible shifts of finite type and Cuntz-Krieger algebras, J. reine angew. Math. **462** (1995), 185–217.
[10] K. H. Kim and F. W. Roush, Some results on decidability of shift equivalence, J. Combin. Inform. System Sci. **4** (1979), 123–146.
[11] ———, Decidability of shift equivalence, in *Dynamical Systems* (J. W. Alexander, ed.), Lecture Notes in Math. **1342** (1988), 374–424, Springer-Verlag, New York.
[12] ———, On the structure of inert automorphisms of subshifts, Pure Math. and Appl. **2** (1991), 3–22.
[13] ———, Williams's conjecture is false for reducible subshifts, JAMS **5** (1992), 213–215.
[14] ———, The Williams conjecture is false for irreducible subshifts, ERA AMS **3** (1997), 105–109.
[15] ———, Topological classification of reducible subshifts, Pure Math. and Appl. **3** (1992), 87–102.
[16] K. H. Kim, F. W. Roush, and J. Wagoner, Automorphisms of the dimension group and gyration numbers, JAMS **5** (1992), 191–212.
[17] ———, Characterization of inert actions on periodic points, preprint (1996).

[18] W. KRIEGER, On dimension functions and topological Markov chains, Invent. Math. **56** (1980), 239–250.
[19] D. LIND and B. MARCUS, *An Introduction to Symbolic Dynamics and Coding*, Cambridge University Press, Cambridge, 1995.
[20] M. NASU, Topological conjugacy for sofic systems and extensions of automorphisms of finite subsystems of topological Markov shifts, Lecture Notes in Math. **1342** (1988), 564–607, Springer-Verlag, New York.
[21] W. PARRY and R. F. WILLIAMS, Block coding and a zeta function for finite Markov chains, Proc. London Math. Soc. **35** (1977), 483–495.
[22] J. B. WAGONER, *Markov Partitions and $K_2$*, Publ. Math. IHES, no. 65 (1987), 91–129.
[23] ______, Triangle identities and symmetries of a subshift of finite type, Pacific J. Math. **144** (1990), 181–205.
[24] ______, Higher-dimensional shift equivalence and strong shift equivalence are the same over the integers, Proc. AMS **109** (1990), 527–536.
[25] ______, Eventual finite order generation for the kernel of the dimension group representation, Trans. AMS **317** (1990), 331–350.
[26] R. F. WILLIAMS, Classification of subshifts of finite type, Ann. of Math. **98** (1973), 120–153; Errata, ibid. **99** (1974), 380–381.
[27] ______, Strong shift equivalence of matrices in GL(2,$\mathbb{Z}$), Contemp. Math. **135** (1992), 445–451.